# Between the LIL and the LSL

ALLAN GUT[1], FREDRIK JONSSON[2] and ULRICH STADTMÜLLER[3]

[1]*Department of Mathematics, Uppsala University, Box 480, SE-751 06 Uppsala, Sweden.*
*E-mail: allan.gut@math.uu.se*

[2]*Department of Mathematics, Uppsala University, Box 480, SE-751 06 Uppsala, Sweden.*
*E-mail: fredrik.jonsson@math.uu.se*

[3]*Ulm University, Department of Number Theory and Probability Theory, D-89069 Ulm, Germany. E-mail: ulrich.stadtmueller@uni-ulm.de*

In two earlier papers, two of the present authors (A.G. and U.S.) extended Lai's [*Ann. Probab.* **2** (1974) 432–440] law of the single logarithm for delayed sums to a multiindex setting in which the edges of the **n**th window grow like $|\mathbf{n}|^\alpha$, or with different $\alpha$'s, where the $\alpha$'s belong to $(0, 1)$. In this paper, the edge of the $n$th window typically grows like $n/\log n$, thus at a higher rate than any power less than one, but not quite at the LIL-rate.

## 1. Introduction

Let $X, \{X_k, k \geq 1\}$ be i.i.d. random variables with mean 0 and partial sums $\{S_n, n \geq 1\}$. The Hartman–Wintner *law of the iterated logarithm* (LIL) states that

$$\limsup_{n \to \infty} (\liminf_{n \to \infty}) \frac{S_n}{\sqrt{2n \log \log n}} = \sigma \ (-\sigma) \qquad \text{a.s.}$$
$$\iff EX^2 < \infty, \qquad EX = 0 \quad \text{and} \quad EX^2 = \sigma^2.$$

The sufficiency was proven by Hartman and Wintner [8], the necessity by Strassen [11].

The *law of the single logarithm* (LSL) is due to Lai [9], and deals with *delayed sums* or *windows*, namely, with

$$T_{n,n+k} = \sum_{j=n+1}^{n+k} X_j, \qquad n \geq 0, k \geq 1,$$

and states that for $0 < \alpha < 1$,

$$\limsup_{n \to \infty} \frac{T_{n,n+n^\alpha}}{\sqrt{2n^\alpha \log n}} = \sigma\sqrt{1-\alpha} \quad \text{a.s.}$$

$$\iff \quad E(|X|^{2/\alpha}(\log^+|X|)^{-1/\alpha}) < \infty, \qquad EX^2 = \sigma^2, \qquad EX = 0,$$

where, throughout, $\log^+ x = \max\{\log x, 1\}$.

The degenerate boundary case $\alpha = 0$ contains the trivial one, in that the window reduces to a single random variable. More precisely, in that case,

$$\frac{T_{n,n+1}}{b_n} = \frac{X_{n+1}}{b_n} \stackrel{a.s.}{\to} 0 \quad \text{as } n \to \infty \iff \sum_{n=1}^{\infty} P(|X| > b_n) < \infty,$$

which, in turn, holds if and only if $Eb^{-1}(|X|) < \infty$, where $b^{-1}(\cdot)$ is a (suitably defined) inverse of $\{b_n\}$.

The next interesting case with $\alpha = 0$ is when the span $a_n = \log n$, that is, the window $T_{n,n+\log n}$, in which case the so-called Erdős–Rényi law ([3], Theorem 2, [2], Theorem 2.4.3) tells us that if $EX = 0$ and the moment generating function $\psi_X(t) = E\exp\{tX\}$ exists in a neigborhood of 0, then for any $c > 0$,

$$\lim_{n \to \infty} \max_{0 \leq k \leq n-k} \frac{T_{k,k+c\log k}}{c\log k} = \rho(c) \qquad \text{a.s.},$$

where

$$\rho(c) = \sup\Big\{x : \inf_t e^{-tx}\psi_X(t) \geq e^{-1/c}\Big\}.$$

Note that here the limit actually depends on the distribution of the summands.

For a generalization to more general window widths $a_n$ such that $a_n/\log n \to \infty$ as $n \to \infty$, but still assuming that the moment generating function exists, see, for example, [2], Theorem 3.1.1, where the limit, in contrast to the result just cited, does not depend on the distribution. Results where the moment condition is somewhat weaker than existence of a moment generating function were discussed in [10]; here, the limit depends on both the variance and the distribution. Using strong invariance principles Lai's result above can be generalized somewhat, see, for example, [2], Theorem 3.2.1, but there at least the $p$th moment, $p > 2$, is needed.

For the boundary case at the other end with $\alpha = 1$, one has $a_n = n$ and $T_{n,2n} \stackrel{d}{=} S_n$, and the correct norming is as in the LIL.

One interesting remaining case is when the window size is larger than any power less than one and, at the same time, not quite linear. This is the starting point of the present paper. Technically, we wish to examine windows of the form

$$T_{n,n+a_n}, \qquad \text{where } a_n = \frac{n}{L(n)} \tag{1.1}$$

with

$$\text{a differentiable function } L(\cdot) \nearrow \infty \in \mathcal{SV} \quad \text{and} \quad \frac{xL'(x)}{L(x)} \searrow \qquad \text{as } x \to \infty. \tag{1.2}$$

*Notation.* $L \in \mathcal{SV}$ means that $L$ is slowly varying at infinity (see, for example, [1] or [5], Section A.7).

The typical case one should have in mind is $L(n) = \log n$, that is, the window $T_{n,n+n/\log n}$.

**Remark 1.1.** Strictly speaking, we should write $a_n = [n/L(n)], a_n = [n/\log n]$ and so on. However, in order to avoid trivial and boring technicalities, we shall treat such sequences as integer-valued whenever convenient.

In Section 2 we present the setup, the main result and the implications for some typical slowly varying functions, namely $L(x) = (\log x)^p$ for $p > 0$ and iterated logarithms $L(x) = \log_m x$, where $\log_m(x)$ denotes the $m$-times iterated logarithm. For the proof, in Section 3, we first review the exponential inequalities. Section 3.2 then introduces a family of subsequences within which sufficiency of the moment condition is proved in Sections 3.3–3.6. Section 3.7 deals with the same issue for the full sequence, while the question of necessity is dealt with in Section 3.8. Proofs of the corollaries in Section 2 are provided in Section 4, while Section 5 furnishes further examples, including some with more complicated slowly varying parts.

It turns out that the proof of the main result has some ingredients in common with that of the classical LIL, primarily in the sense that one needs two truncations, one to match the Kolmogorov exponential bounds and one to match the moment requirements. Typically (and somewhat frustratingly), it is the thin central part that causes the main trouble in the proof. A weaker result is obtained if only the first truncation is made. The cost is that too much integrability will be required. However, for the reader who is not so concerned with optimality, we include a proof of this weaker version in Section 6, after which we revisit two examples in order to illustrate the consequences.

## 2. Setup and main result

Recall that the window widths, $a_n$, are assumed to be of the form $n/L(n)$, where the function $L$ satisfies (1.2). Define $\{d_n, n \geq 2\}$ by

$$d_n = \log \frac{n}{a_n} + \log \log n = \log L(n) + \log \log n.$$

Note that $\{d_n\}$ may be viewed as the *additional norming sequence* in Theorem 2.1, in the sense that it corresponds to $\{\log_2 n\}$ in the LIL and $\{\log n\}$ in the LSL.

Furthermore, let

$$f(n) = \min\{a_n \cdot d_n, n\},$$

with $f(\cdot)$ an increasing interpolating function, that is, $f(x) = f_{[x]}$ for $x > 0$ and $f^{-1}(\cdot)$, the corresponding (suitably defined) inverse function.

Here, now, is our main result.

**Theorem 2.1.** *Suppose that $X, X_1, X_2, \ldots$ are i.i.d. random variables with mean $0$ and finite variance $\sigma^2$, and let $T_{n,n+k} = \sum_{j=n+1}^{n+k} X_j$. If*

$$E(f^{-1}(X^2)) < \infty, \tag{2.1}$$

*then*

$$\limsup_{n \to \infty} \frac{T_{n,n+a_n}}{\sqrt{2a_n d_n}} = \sigma \qquad a.s. \tag{2.2}$$

*Conversely, if*

$$P\left(\limsup_{n \to \infty} \frac{|T_{n,n+a_n}|}{\sqrt{a_n d_n}} < \infty\right) > 0, \tag{2.3}$$

*then (2.1) holds, $EX = 0$ and (2.2) holds with $\sigma^2 = \operatorname{Var} X$.*

*Remark 2.1.* The "natural" necessary moment assumption is (2.1) with $f(n) = a_n d_n$. However, for very slowly increasing functions $L$, for example, $L(x) = \log\log\log\log x$, it turns out that finite variance is needed, and since we then have $f(n) = n$, (2.1) is equivalent to finite variance.

*Remark 2.2.* The result also holds for any sequence $\{a_n\}$ which is of regular variation of order $\alpha \in (0, 1)$. Here, the sufficiency part can be obtained from strong invariance principles, as described in, for example, Theorem 3.2.2 in the book [2]. However, for our situation, no strong invariance principle is available.

*Remark 2.3.* In addition to the $\limsup$ results, there exist, throughout, $\liminf$ counterparts such that $\liminf \cdots = -\limsup \cdots$ a.s. Actually, the set of limit points is the whole interval $[-\sigma, \sigma]$.

The slowly varying function that immediately comes to mind is (of course) the logarithmic function. The second one would be $L(x) = \log_2 x = \log\log x$ and, possibly, $L(x) = \log_m(x)$. We precede the proofs by stating the conclusions for these cases as separate corollaries. For simplicity, we omit the converse parts.

**Corollary 2.1.** *Suppose that $X, X_1, X_2, \ldots$ are i.i.d. random variables with mean $0$ and finite variance $\sigma^2$, and let $T_{n,n+k} = \sum_{j=n+1}^{n+k} X_j$. If, for some $p > 0$,*

$$EX^2 \frac{(\log^+ |X|)^p}{\log^+ \log^+ |X|} < \infty, \tag{2.4}$$

*then*

$$\limsup_{n \to \infty} \frac{T_{n,n+n/(\log n)^p}}{\sqrt{2(p+1)(n/(\log n)^p)\log\log n}} = \sigma \qquad a.s. \tag{2.5}$$

**Corollary 2.2.** *Suppose that $X, X_1, X_2, \ldots$ are i.i.d. random variables with mean 0 and let $T_{n,n+k} = \sum_{j=n+1}^{n+k} X_j$. If $\sigma^2 = \operatorname{Var} X < \infty$, then for any $m \geq 2$,*

$$\limsup_{n \to \infty} \frac{T_{n,n+n/\log_m(n)}}{\sqrt{2(n/\log_m(n)) \log \log n}} = \sigma \qquad a.s. \tag{2.6}$$

Note that in the case $m = 2$, the normalization is just $\sqrt{2n}$. Proofs of the corollaries are deferred to Section 4 and further examples are given in Section 5.

## 3. Proof of Theorem 2.1

In spite of the fact that we are dealing with limit laws for delayed sums, the present topic is, in fact, too close to the LIL to warrant LSL techniques. In contrast to the proofs in [9] and [6, 7], where one uses exponential bounds and Borel–Cantelli lemmas for the single primed contribution along a suitably subsequence, and takes care of the double and triple primed contributions for the full sequence and fills the gaps, we have to resort to the LIL technique where one proves Borel–Cantelli lemmas and thus also the theorem itself, first for subsequences and then for the entire sequence.

We thus begin by providing Borel–Cantelli sums along subsequences, after which an appeal to the Borel–Cantelli lemmas completes the proof for subsequences.

Section 3.7 is devoted to the problem of "filling the gaps" in order to include arbitrary windows.

### 3.1. Truncation and exponential bounds

The typical approach to proving results of the LIL type requires two truncations: the first to match the Kolmogorov exponential bounds (see, for example, [5], Section 8.2) and the second to match the moment requirements.

To this end, we introduce parameters $\delta > 0$ and $\varepsilon > 0$, and let

$$b_n = \frac{\sigma \delta}{\varepsilon} \sqrt{\frac{a_n}{d_n}} \tag{3.1}$$

and

$$X'_n = X_n I\{|X_n| \leq b_n\}, \qquad X''_n = X_n I\{b_n < |X_n| < \delta \sqrt{f(n)}\},$$
$$X'''_n = X_n I\{|X_n| \geq \delta \sqrt{f(n)}\}.$$

In the following, all objects with primes or multiple primes refer to the respective truncated summands.

Since truncation destroys centering, we obtain, using standard procedures, together with the fact that $EX = 0$,

$$|EX'_k| = |-EX_k I\{|X_k| > b_k\}| \leq E|X|I\{|X_k| > b_k\} \leq \frac{EX^2 I\{|X| > b_k\}}{b_k}$$

so that

$$\begin{aligned}
|ET'_{n,n+a_n}| &\leq \sum_{n \leq k \leq n+a_n} \frac{EX^2 I\{|X| > b_k\}}{b_k} \leq a_n \cdot \frac{EX^2 I\{|X| > b_n\}}{b_n} \\
&= \frac{\varepsilon}{\sigma\delta} \cdot \sqrt{a_n d_n} \cdot EX^2 I\{|X| > b_n\} = \mathrm{o}(\sqrt{a_n d_n}) \qquad \text{as } n \to \infty.
\end{aligned} \tag{3.2}$$

*Upper bounds*

Since

$$\operatorname{Var} X'_k \leq E(X'_k)^2 \leq EX^2 = \sigma^2,$$

it follows that

$$\operatorname{Var}(T'_{n,n+a_n}) \leq a_n \sigma^2. \tag{3.3}$$

An application of the Kolmogorov upper exponential bound (see, for example, [5], Lemma 8.2.1) with $x = \varepsilon(1-\delta)\sqrt{2d_n}$ and $c_n = 2\delta/x$, together with (3.2) and (3.3), now yields

$$\begin{aligned}
P(T'_{n,n+a_n} > \varepsilon\sqrt{2a_n d_n}) &\leq P(T'_{n,n+a_n} - ET'_{n,n+a_n} > \varepsilon(1-\delta)\sqrt{2a_n d_n}) \\
&\leq P\left(T'_{n,n+a_n} - ET'_{n,n+a_n} > \frac{\varepsilon(1-\delta)}{\sigma}\sqrt{2\operatorname{Var}(T'_{n,n+a_n})d_n}\right) \\
&\leq \exp\left\{-\frac{2\varepsilon^2(1-\delta)^2}{2\sigma^2} \cdot d_n(1-\delta)\right\} \\
&= \exp\left\{-\frac{\varepsilon^2(1-\delta)^3}{\sigma^2} \cdot d_n\right\}.
\end{aligned} \tag{3.4}$$

*Lower bounds*

In order to apply the lower exponential bound (see, for example, [5], Lemma 8.2.2), we first need a lower bound for the truncated variances:

$$\begin{aligned}
\operatorname{Var} X'_k &= E{X'_k}^2 - (EX'_k)^2 = EX^2 - EX^2 I\{|X_k| \geq b_k\} - (EX'_k)^2 \\
&\geq \sigma^2 - 2EX^2 I\{|X_k| \geq b_k\} \geq \sigma^2(1-\delta)
\end{aligned}$$

for $n$ large, so that

$$\operatorname{Var}(T'_{n,n+a_n}) \geq a_n \sigma^2(1-\delta) \qquad \text{for } n \text{ large.} \tag{3.5}$$

It now follows that for any $\gamma > 0$,

$$
\begin{aligned}
P(T'_{n,n+a_n} > \varepsilon\sqrt{2a_n d_n}) &\geq P(T'_{n,n+a_n} - ET'_{n,n+a_n} > \varepsilon(1+\delta)\sqrt{2a_n d_n}) \\
&\geq P\left(T'_{n,n+a_n} - ET'_{n,n+a_n} > \frac{\varepsilon(1+\delta)}{\sigma\sqrt{(1-\delta)}}\sqrt{2\operatorname{Var}(T'_{n,n+a_n})d_n}\right) \\
&\geq \exp\left\{-\frac{2\varepsilon^2(1+\delta)^2}{2\sigma^2(1-\delta)}\cdot d_n(1+\gamma)\right\} \\
&= \exp\left\{-\frac{\varepsilon^2(1+\delta)^2(1+\gamma)}{\sigma^2(1-\delta)}\cdot d_n\right\} \qquad \text{for } n \text{ large.}
\end{aligned}
\tag{3.6}
$$

### 3.2. A family of subsequences

In order to choose a suitable subsequence, consider the difference equation $n_{k+1} - n_k = cn_k/L(n_k)$ with a suitable constant $c > 0$ to be determined later, or, in continuous variables,

$$y' = cy/L(y). \tag{3.7}$$

With $\varphi(y) = \int^y \frac{L(u)\,du}{u}$ being in the class $\Pi$ (see [1] for the notation and Theorem 3.7.3) and $\psi(x) = \varphi^{-1}(x)$ being in the class $\Gamma$ (see [1] for the notation and Theorem 3.10.4), the solution of the differential equation is given by $\psi(cx)$ and the subsequence of interest is $n_k = \psi(ck)$. Note that $\frac{n_{k+1}}{n_k} = 1 + \frac{c}{L(n_k)} \to 1$ and that $L(n_{k+1})/L(n_k) \to 1$ as $k \to \infty$.

An important relation in the following is

$$d_{n_k} = \log(L(\psi(ck))\log\psi(ck)) \sim \log(ck) \sim \log k \qquad \text{as } t \to \infty \text{ for any } c > 0, \tag{3.8}$$

which is an immediate consequence of the following result.

**Lemma 3.1.** *With a slowly varying function $L(\cdot)$ satisfying (1.2), we have*

$$\frac{\log(L(t)\log t)}{\log\varphi(t)} \to 1 \qquad \text{as } t \to \infty. \tag{3.9}$$

**Proof.** With $\varphi^*(t) = L(t)\log t$, we have $\varphi(t) \leq \varphi^*(t)$ since $L(\cdot) \nearrow$. Next,

$$
\begin{aligned}
\varphi^*(t) &= \int_1^t \left(L'(u)\log u + \frac{L(u)}{u}\right)du = \int_1^t \frac{L'(u)uL(u)}{L(u)u}\int_1^u \frac{1}{v}\,dv\,du + \varphi(t) \\
&\leq \int_1^t \frac{L(u)}{u}\int_1^u \frac{L'(v)}{L(v)}\,dv\,du + \varphi(t) \leq \varphi(t)(1 + \log(L(t))),
\end{aligned}
$$

where we used condition (1.2). Hence,

$$1 \geq \frac{\log\varphi(t)}{\log\varphi^*(t)} \geq 1 - \frac{\log(1+\log L(t))}{\log(L(t)\log t)} \to 1 \qquad \text{as } t \to \infty. \qquad \square$$

***Remark 3.1.*** In the classical proof of the LIL, the subsequence has a geometric growth rate: $\{\lambda^k, k \geq 1\}$ for some $\lambda$ close to 1. For the LSL, the subsequence has a polynomial growth rate: $\{(k/\log k)^{1/(1-\alpha)}, k \geq 1\}$. It is therefore natural in the present, intermediate, context to search for a subsequence with a growth rate between geometric and polynomial, that is, to search for something like $\{\lambda^{k^\beta}, k \geq 1\}$ for some $\beta \in (0,1)$. For the canonical case $L(n) = \log n$, it turns out that $n_k \sim e^{c\sqrt{2k}}$.

## 3.3. Sufficiency along subsequences: $T'_{n,n+a_n}$

*The upper bound*

Here, we use $c > 0$ small. Let $\{n_k = \psi(ck), k \geq 1\}$, where $n_k \nearrow \infty$ as $k \to \infty$, satisfy

$$\sum_{k=1}^{\infty} \exp\left\{-\frac{\varepsilon^2(1-\delta)^3}{\sigma^2} \cdot d_{n_k}\right\} < \infty. \tag{3.10}$$

Applying (3.4) to $\{X'_k, k \geq 1\}$ then yields

$$\sum_{k=1}^{\infty} P(|T'_{n_k, n_k+a_{n_k}}| > \varepsilon\sqrt{2a_{n_k}d_{n_k}}) < \infty \tag{3.11}$$

for any $\varepsilon > \sigma$. Note that (3.11) is independent of the special choice of $c > 0$.

*The lower bound*

We now choose the sparser subsequence $\{n_k = \psi(ck), k \geq 1\}$, where $c > 1$ and $n_k \nearrow \infty$ as $k \to \infty$, satisfying

$$\sum_{k=1}^{\infty} \exp\left\{-\frac{\varepsilon^2(1+\delta)^2(1+\gamma)}{\sigma^2(1-\delta)} \cdot d_{n_k}\right\} = \infty \tag{3.12}$$

for any $\varepsilon < \sigma$. Observe that the windows are now non-overlapping since $c > 1$ implies that $n_{k+1} > n_k + n_k/L(n_k)$ eventually. Applying (3.6) to this sequence similarly shows that

$$\sum_{k=1}^{\infty} P(T'_{n_k, n_k+a_{n_k}} > \varepsilon\sqrt{2a_{n_k}d_{n_k}}) = \infty. \tag{3.13}$$

## 3.4. Sufficiency along subsequences: $T''_{n,n+a_n}$

The next step is to prove the analog of (3.11) for $T''_{n,n+a_n}$, that is,

$$\sum_{k=1}^{\infty} P(|T''_{n,n+a_n}| > \delta\sqrt{f(n)}) < \infty. \tag{3.14}$$

*The symmetric case*

We first consider symmetric random variables, beginning by recalling the Kahane–Hoffmann–Jørgensen inequality (see, for example, [5], Theorem 3.7.5).

**Lemma 3.2.** *Suppose that $X_1, X_2, \ldots, X_n$ are independent symmetric random variables with partial sums $S_n$, $n \geq 1$.*

(i) *For any $x, y > 0$,*

$$P(|S_n| > 2x + y) \leq P\Big(\max_{1 \leq k \leq n} |X_k| > y\Big) + 4(P(|S_n| > x))^2$$

$$\leq \sum_{k=1}^{n} P(|X_k| > y) + 4(P(|S_n| > x))^2.$$

(ii) *If $X_1, X_2, \ldots, X_n$ are identically distributed (and $x = y$), then an iteration yields that there are constants $\kappa_i > 0$, $i = 1, 2$, such that*

$$P(|S_n| > 9x) \leq \kappa_1 n P(|X_1| > x) + \kappa_2 (P(|S_n| > x))^4.$$

Applying the lemma to $T''_{n,n+a_n}$, we thus obtain, with $\eta = \delta/9$, that

$$P(|T''_{n,n+a_n}| > \delta\sqrt{f(n)})$$
$$\leq \kappa_1 \sum_{k=n+1}^{n+a_n} P(|X''_k| > \eta\sqrt{f(n)}) + \kappa_2 (P(|T''_{n,n+a_n}| > \eta\sqrt{f(n)}))^4 \qquad (3.15)$$
$$\leq \kappa_1 a_n P(|X| > \eta\sqrt{f(n)}) + \kappa_2 (P(|T''_{n,n+a_n}| > \eta\sqrt{f(n)}))^4.$$

Summing over our subsequence for $k_0$ large (remembering that $d_{n_k} \sim \log k$ as $k \to \infty$), we now have

$$\sum_{k=k_0}^{\infty} a_{n_k} P(|X| > \eta\sqrt{f(n_k)})$$
$$\leq \sum_{k=k_0}^{\infty} a_{n_k} P(f^{(-1)}(X^2/\eta^2) > n_k) = \sum_{k=k_0}^{\infty} \frac{n_k}{L(n_k)} P(f^{(-1)}(X^2/\eta^2) > n_k)$$
$$\leq \int_1^{\infty} \frac{\psi(x)}{L(\psi(x))} P(f^{(-1)}(X^2/\eta^2) > c\psi(x)) \, dx$$
$$\qquad \qquad \qquad \qquad \qquad \qquad \qquad \qquad \qquad \qquad \qquad (3.16)$$
$$\left(\text{use } \frac{\psi(x-1)}{\psi(x)} \geq c > 0, \text{ a change of variable } y = \psi(x), \text{ hence, } \frac{dx}{dy} = \varphi'(y) = \frac{L(y)}{y}\right)$$
$$= \int_C^{\infty} \frac{y}{L(y)} P(f^{(-1)}(X^2/\eta^2) > cy) \frac{L(y)}{y} \, dy = \int_C^{\infty} P(f^{(-1)}(X^2/\eta^2) > cy) \, dy$$

$$\leq C E f^{(-1)}(X^2) < \infty,$$

which takes care of the first term in the right-hand side of (3.15).

As for the second one, Chebyshev's inequality tells us that

$$P(|T''_{n,n+a_n}| > \eta\sqrt{f(n)}) \leq \frac{\operatorname{Var} T''_{n,n+a_n}}{\eta^2 f(n)} \leq \frac{a_n E X^2 I\{b_n < |X| < \delta\sqrt{f(n+n/L(n))}\}}{\eta^2 f(n)}$$

$$= \frac{\varepsilon}{\sigma\delta} \frac{E X^2 I\{b_n < |X| < \delta\sqrt{f(n+n/L(n))}\}}{\eta^2 d_n} \leq \frac{\varepsilon}{\sigma\delta} \frac{E X^2}{\eta^2 d_n}$$

and, hence, that

$$(P(|T''_{n,n+a_n}| > \eta\sqrt{f(n)}))^4 \leq \left(\frac{\varepsilon}{\sigma\delta\eta^2 d_n}\right)^4 (E X^2)^3 E X^2 I\{b_n < |X| < \delta\sqrt{f(n+n/L(n))}\}$$

so that

$$\left(\frac{\sigma\delta\eta^2}{\varepsilon}\right)^4 \frac{1}{(EX^2)^3} \sum_{k=k_0}^{\infty} P(|T''_{n_k,n_k+n_k/L(n_k)}| > \eta\sqrt{f(n_k)}) \quad (3.17)$$

$$\leq C \sum_{k=k_0}^{\infty} \frac{E X^2 I\{b_{n_k} < |X_k| < \delta\sqrt{f(n_k+n_k/L(n_k))}\}}{d_{n_k}^4}$$

$$\leq C \sum_{k=k_0}^{\infty} \frac{1}{(\log k)^4} \int_{b_{n_k}}^{\delta\sqrt{f(n_k+n_k/L(n_k))}} x^2 \, dF(x)$$

$$= \int_{k_*}^{\infty} \left(\sum_{A(k,x)} \frac{1}{(\log k)^4}\right) x^2 \, dF(x), \quad (3.18)$$

where $k_*$ is some irrelevant lower limit and

$$A(k,x) = \{k : b_{n_k} < |x| < \delta\sqrt{f(n_k+n_k/L(n_k))}\}.$$

In order to invert the double inequality, we first observe that in the case $f(n) = a_n d_n$ (the case $f(n) = n$ is simpler and only the necessary changes are indicated),

$$a_{n_k+n_k/L(n_k)} = \frac{n_k + n_k/L(n_k)}{L(n_k + n_k/L(n_k))} \leq \frac{n_k}{L(n_k)} \left(1 + \frac{1}{L(n_k)}\right)$$

and that

$$d_{n_k+n_k/L(n_k)} = \log L(n_k+n_k/L(n_k)) + \log\log(n_k+n_k/L(n_k)) \sim \log d_{n_k} \sim \log k \sim \log(\varphi(n_k))$$

because of the slow variation of $L$, $\log L$, and $\log x$, and the fact that we have chosen our subsequence via the relation $n_k = \psi(ck)$, which implies that $\varphi(n_k) \sim ck$.

Exploiting this yields

$$f(n_k + n_k/L(n_k)) = a_{n_k+n_k/L(n_k)} \cdot d_{n_k+n_k/L(n_k)}$$
$$\leq (1+\delta/2)\frac{n_k}{L(n_k)}\left(1+\frac{1}{L(n_k)}\right)\cdot \log(\varphi(n_k)) \quad (3.19)$$
$$\leq \frac{n_k}{L(n_k)}(1+\delta)\cdot \log(\varphi(n_k)) \quad \text{as } k\to\infty.$$

Next, for a slowly varying function $L$, let $L^{\#}$ be its de Bruijn conjugate (see, for example, [1], Section 1.5), obeying

$$L(xL^{\#}(x))L^{\#}(x)\to 1 \quad \text{and} \quad L(x)L^{\#}(xL(x))\to 1 \quad \text{as } x\to\infty. \quad (3.20)$$

With its help, we can solve $t=\xi L(\xi)$ asymptotically by $\xi\sim tL^{\#}(t)$. Now, for evaluating $A(k,x)$, we define

$$L_1(u)=\frac{1}{L(u)\log(\varphi(u))} \quad \text{and} \quad L_2(u)=\frac{\log(\varphi(u))}{L(u)},$$

both of which are slowly varying (in the case $f(n)=n$, we may define $L_2\equiv 1$), and their de Bruijn conjugates $L_1^{\#}(x)$ and $L_2^{\#}(x)$. Then, with suitable constants $c_i$, we have

$$A(k,x)\subset\left\{k:\left(\frac{\delta\sigma}{\varepsilon}\right)^2\frac{n_k}{L(n_k)\log(\varphi(n_k))}\leq x^2\leq \delta^2\frac{n_k}{L(n_k)}(1+\delta)\cdot\log(\varphi(n_k))\right\}$$
$$\subset\left\{k:\left(\frac{\delta\sigma}{\varepsilon}\right)^2\frac{\psi(ck)}{L(\psi(ck))\log(\varphi(\psi(ck)))}\leq x^2\leq \delta^2\frac{\psi(ck)}{L(\psi(ck))}(1+\delta)\cdot\log(\varphi(\psi(ck)))\right\}$$
$$\subset\{k:c_1\psi(ck)L_1(\psi(ck))\leq x^2\leq c_2\psi(ck)L_2(\psi(ck))\}$$
$$\subset\{k:c_3x^2L_2^{\#}(x^2)\leq \psi(ck)=n_k\leq c_4x^2L_1^{\#}(x^2)\}$$
$$\subset\{k:c_5\varphi(x^2L_2^{\#}(x^2))\leq k\leq c_6\varphi(x^2L_1^{\#}(x^2))\}.$$

An application of the mean value theorem, the fact that $\varphi'(x)=L(x)/x\searrow$ as $x\to\infty$ and (3.20) therefore imply that

$$\mathrm{Card}(A(k,x))\leq c_6\varphi(x^2L_1^{\#}(x^2))-c_5\varphi(x^2L_2^{\#}(x^2))\leq C(\varphi(x^2L_1^{\#}(x^2))-\varphi(x^2L_2^{\#}(x^2)))$$
$$\leq C\varphi'(x^2L_2^{\#}(x^2))(x^2L_1^{\#}(x^2)-x^2L_2^{\#}(x^2))$$
$$=C\frac{L(x^2L_2^{\#}(x^2))}{x^2L_2^{\#}(x^2)}(x^2L_1^{\#}(x^2)-x^2L_2^{\#}(x^2))$$
$$\leq C\frac{\log(\varphi(x^2L_2^{\#}(x^2)))}{L_2(x^2L_2^{\#}(x^2))L_2^{\#}(x^2)}(L_1^{\#}(x^2)-L_2^{\#}(x^2))$$
$$\leq C\log(\varphi(x^2L_2^{\#}(x^2)))L_1^{\#}(x^2).$$

Inserting this into the inner sum in (3.18), we now obtain

$$\sum_{A(k,x)} \frac{1}{(\log k)^4} \leq C \frac{L_1^{\#}(x^2)}{(\log \varphi(x^2 L_1^{\#}(x^2)))^3} \leq C L_2^{\#}(x^2), \tag{3.21}$$

where the last inequality follows from the fact that $\frac{L_1^{\#}(x^2)}{(\log \varphi(x^2 L_1^{\#}(x^2)))^3} \leq C L_2^{\#}(x^2)$ since, using (3.20),

$$\begin{aligned}
\frac{L_1^{\#}(x)}{L_2^{\#}(x)} &\sim \frac{L_2(x L_2^{\#}(x))}{L_1(x L_1^{\#}(x))} \\
&\leq \frac{L(x L_1^{\#}(x))}{L(x L_2^{\#}(x))} \cdot (\log(\varphi(x L_1^{\#}(x))))^2 \\
&\leq C (\log(\varphi(x L_1^{\#}(x))))^2 \exp\left( \int_{x L_2^{\#}(x)}^{x L_1^{\#}(x)} \varepsilon(t)/t \, dt \right) \\
&\leq C (\log(\varphi(x L_1^{\#}(x))))^2 \exp\left( o(1) \log\left( \frac{L_1^{\#}(x)}{L_2^{\#}(x)} \right) \right),
\end{aligned}$$

by the representation theorem for slowly varying functions. (In the case $f(n) = n$, the inequality (3.21) is trivial since $L_2^{\#} \equiv 1$ and $L_1^{\#}$ is decreasing.) Finally, using the fact that $f^{-1}(x) \sim x L_2^{\#}(x)$, we conclude that the sum in (3.18) converges, which takes care of the second sum in (3.15).

Combining this with (3.16) proves the validity of (3.14) in the symmetric case.

*Desymmetrization*

In order to prove (3.14) for the general case, we first estimate the truncated means. Remembering that $EX_k = 0$ for all $k$, we obtain, by stretching the bounds to the extreme,

$$\begin{aligned}
|ET''_{n,n+a_n}| &= \left| \sum_{k=n+1}^{n+a_n} EX_k I\{b_k < |X_k| < \delta \sqrt{f_k}\} \right| \\
&\leq \sum_{k=n+1}^{n+a_n} E|X_k| I\{b_n < |X_k| < \delta \sqrt{f(n+n/L(n))}\} \\
&\leq a_n E|X| I\{|X| \geq b_n\} \leq \frac{a_n}{b_n} EX^2 I\{|X| \geq b_n\} \\
&= \frac{\varepsilon}{\sigma \delta} \sqrt{a_n d_n} EX^2 I\{|X| \geq b_n\} = o(\sqrt{a_n d_n}) \qquad \text{as } n \to \infty
\end{aligned}$$

since this is the same estimate as for $ET'_{n,n+a_n}$, after which the desired conclusion follows with the aid of the symmetrization inequalities (see [5], Proposition 3.6.2).

## 3.5. Sufficiency along subsequences: $T'''_{n,n+a_n}$

In order for $|T'''_{n,n+a_n}|$ to surpass the level $\eta\sqrt{a_n d_n}$, it is necessary that at least one of the $X'''$'s is non-zero. For every $\eta > 0$ (recall that $a_{n_k} = n_k/L(n_k)$, $d_{n_k} \sim \log k$), this means that

$$\sum_{k=1}^{\infty} P(|T'''_{n_k,n_k+n_k/L(n_k)}| > \eta\sqrt{a_{n_k} d_{n_k}}) \leq \sum_{k=1}^{\infty}\sum_{j=1}^{a_{n_k}} P\left(|X_{n_k+j}| > \frac{\eta}{2}\sqrt{f(n_k+j)}\right)$$
$$\leq \sum_{k=1}^{\infty} a_{n_k} P\left(|X| > \frac{\eta}{2}\sqrt{f(n_k)}\right) < \infty, \tag{3.22}$$

by (3.16).

## 3.6. Sufficiency along subsequences: Combining the contributions

Combining (3.11), (3.14) and (3.22), we conclude that

$$\sum_{k=1}^{\infty} P(|T_{n_k,n_k+n_k/L(n_k)}| > (\varepsilon+2\eta)\sqrt{2a_{n_k} d_{n_k}}) < \infty \tag{3.23}$$

provided $\varepsilon > \sigma/(1-\delta)^{3/2}$ and, since $\eta$ and $\delta$ may be arbitrarily chosen, that

$$\sum_{k=1}^{\infty} P(|T_{n_k,n_k+n_k/L(n_k)}| > \varepsilon\sqrt{2a_{n_k} d_{n_k}}) < \infty \qquad \text{for } \varepsilon > \sigma \tag{3.24}$$

so that, in view of the first Borel–Cantelli lemma,

$$\limsup_{k \to \infty} \frac{T_{n_k,n_k+n_k/L(n_k)}}{\sqrt{2a_{n_k} d_{n_k}}} \leq \sigma \qquad \text{a.s.} \tag{3.25}$$

A completely analogous argument, combining (3.13), (3.14) and (3.22), yields

$$\sum_{k=1}^{\infty} P(T_{n_k,n_k+n_k/L(n_k)} > \varepsilon\sqrt{2a_{n_k} d_{n_k}}) = \infty \qquad \text{for } \varepsilon < \sigma \tag{3.26}$$

and since the windows with this, sparser, subsequence are disjoint, we may apply the second Borel–Cantelli lemma to conclude that

$$\limsup_{k \to \infty} \frac{T_{n_k,n_k+n_k/L(n_k)}}{\sqrt{2a_{n_k} d_{n_k}}} \geq \sigma \qquad \text{a.s.} \tag{3.27}$$

Finally, combining (3.25) and (3.27) yields

$$\limsup_{k\to\infty} \frac{T_{n_k, n_k+n_k/L(n_k)}}{\sqrt{2a_{n_k}d_{n_k}}} = \sigma \qquad \text{a.s.}, \tag{3.28}$$

which, in addition, proves the sufficiency of the following result, which is Theorem 2.1 for subsequences of the form $n_k = \psi(ck)$ with $c > 1$. The necessity follows, of course, from the necessity for the full sequence, the proof of which is given in Section 3.8 below.

**Theorem 3.1.** *Suppose that $X, X_1, X_2, \ldots$ are i.i.d. random variables with mean $0$ and finite variance $\sigma^2$, let $T_n = \sum_{j=n+1}^{n+k} X_j$ and let, for $c > 1$,*

$$n_k = \varphi^{(-1)}(ck), \qquad k \geq 1,$$

*where $\varphi^{(-1)}$ is the inverse of $\varphi(y) = \int^y \frac{L(u)}{u}\,du$. If (2.1) holds, then*

$$\limsup_{n\to\infty} \frac{T_{n_k + n_k/L(n_k)}}{\sqrt{2(n_k/L(n_k))\log k}} = \sigma \qquad a.s. \tag{3.29}$$

*Conversely, if*

$$P\left(\limsup_{n\to\infty} \frac{|T_{n_k + n_k/L(n_k)}|}{\sqrt{(n_k/L(n_k))\log k}} < \infty\right) > 0, \tag{3.30}$$

*then (2.1) holds, $EX = 0$, $EX^2 < \infty$ and (3.29) holds with $\sigma^2 = \operatorname{Var} X$.*

### 3.7. Sufficiency for the entire sequence

We must thus show that our process behaves accordingly for the entire sequence. Here, the second Lévy inequality (see, for example, [5], Theorem 3.7.2) is instrumental. Let $\eta > 0$ be given. Then,

$$\begin{aligned}
P&\left(\max_{n_k \leq n \leq n_{k+1}} \frac{S_{n+a_n} - S_n}{\sqrt{2a_n d_n}} > (1+6\eta)\sigma\right) \\
&\leq P\left(\max_{n_k \leq n \leq n_{k+1}} (S_{n+a_n} - S_{n_k+a_{n_k}}) > 2\eta\sigma\sqrt{2a_{n_k}d_{n_k}}\right) \\
&\quad + P\left(\max_{n_k \leq n \leq n_{k+1}} (-S_n + S_{n_k}) > 2\eta\sigma\sqrt{2a_{n_k}d_{n_k}}\right) \\
&\quad + P\left(\max_{n_k \leq n \leq n_{k+1}} (S_{n_k+a_{n_k}} - S_{n_k}) > (1+2\eta)\sigma\sqrt{2a_{n_k}d_{n_k}}\right),
\end{aligned} \tag{3.31}$$

where $n_k = \psi(cn_k)$ as before with some suitable constant $c > 0$ to be fixed shortly.

*LIL and LSL* 15

Set $\tilde{n}_k = n_k + a_{n_k}$. Since $n_{k+1}/n_k \to 1$ and $L(\cdot) \in \mathcal{SV}$, the following relations hold eventually (that is, for $k$ sufficiently large):

$$n_{k+1} - n_k \leq c\psi'(c(k+1)) = c\frac{n_{k+1}}{L(cn_k)} \leq 2ca_{n_k},$$

$$n_k \leq \tilde{n}_k = n_k(1 + (L(n_k))^{-1}) \leq n_k(1+\eta),$$

$$a_{n_k} \leq a_{\tilde{n}_k} \leq a_{n_k(1+\eta)} \leq (1+\eta)a_{n_k},$$

$$\tilde{n}_{k+1} - \tilde{n}_k \leq (n_{k+1} - n_k)(1 + (L(n_k))^{-1}) \leq 2c(1+\eta)a_{n_k} \leq 2c(1+\eta)a_{\tilde{n}_k},$$

$$d_{n_k} \leq d_{\tilde{n}_k} \leq (1+\eta)d_{n_k}.$$

In the following, we exploit these relations without specifically mentioning them each time.

As a first application, we note that (3.31) can be bounded by

$$\leq P\Big(\max_{\tilde{n}_k \leq n \leq \tilde{n}_k + 2c(1+\eta)a_{\tilde{n}_k}} (S_n - S_{\tilde{n}_k}) > 2\eta\sigma\sqrt{2a_{n_k}d_{n_k}}\Big)$$
$$+ P\Big(\max_{n_k \leq n \leq n_k + 2ca_{n_k}} (-S_n + S_{n_k}) > 2\eta\sigma\sqrt{2a_{n_k}d_{n_k}}\Big) \qquad (3.32)$$
$$+ P\Big(\max_{n_k \leq n \leq n_{k+1}} (S_{n_k+a_{n_k}} - S_{n_k}) > (1+2\eta)\sigma\sqrt{2a_{n_k}d_{n_k}}\Big).$$

Now,

$$\mathrm{Var}(S_{\tilde{n}_k+2c(1+\eta)a_{\tilde{n}_k}} - S_{\tilde{n}_k}) = 2c(1+\eta)a_{\tilde{n}_k}\sigma^2 = \mathrm{o}(a_{n_k}d_{n_k}) \qquad \text{as } k \to \infty,$$

$$\mathrm{Var}(S_{n_k+2ca_{n_k}} - S_{n_k}) = 2ca_{n_k}\sigma^2 = \mathrm{o}(a_{n_k}d_{n_k}) \qquad \text{as } k \to \infty,$$

$$\mathrm{Var}(S_{n_k+a_{n_k}} - S_{n_k}) = a_{n_k}\sigma^2 = \mathrm{o}(a_{n_k}d_{n_k}) \qquad \text{as } k \to \infty,$$

that is, the variances are $\leq \eta^4\sigma^2 a_{n_k}d_{n_k}$ for $k$ sufficiently large.

An application of the Lévy inequality to the first two probabilities in (3.32), leaving the third one as is, then shows that (3.32) can be bounded by

$$\leq 2P((S_{\tilde{n}_k+2c(1+\eta)a_{\tilde{n}_k}} - S_{\tilde{n}_k}) > 2\eta\sigma\sqrt{2a_{n_k}d_{n_k}} - \sqrt{2}\cdot\eta^2\sigma\sqrt{a_{n_k}d_{n_k}})$$
$$+ 2P(-(S_{n_k+2ca_{n_k}} - S_{n_k}) > 2\eta\sigma\sqrt{2a_{n_k}d_{n_k}} - \sqrt{2}\cdot\eta^2\sigma\sqrt{a_{n_k}d_{n_k}})$$
$$+ 2P((S_{n_k+a_{n_k}} - S_{n_k}) > (1+\eta)\sigma\sqrt{2a_{n_k}d_{n_k}})$$
$$\leq 2P((S_{\tilde{n}_k+2c(1+\eta)a_{\tilde{n}_k}} - S_{\tilde{n}_k}) > \eta\sigma\sqrt{2a_{n_k}d_{n_k}})$$
$$+ 2P(-(S_{n_k+2c(1+\eta)a_{n_k}} - S_{n_k}) > \eta\sigma\sqrt{2a_{n_k}d_{n_k}}) \qquad (3.33)$$
$$+ 2P((S_{n_k+a_{n_k}} - S_{n_k}) > (1+\eta)\sigma\sqrt{2a_{n_k}d_{n_k}})$$

$$\leq 2P\bigg(|S_{\tilde{n}_k+2c(1+\eta)a_{\tilde{n}_k}} - S_{\tilde{n}_k}| > \frac{\eta}{\sqrt{2c(1+\eta)^3}}\sigma\sqrt{2 \cdot 2c(1+\eta)a_{\tilde{n}_k}d_{\tilde{n}_k}}\bigg)$$

$$+ 2P\bigg(|S_{n_k+2ca_{n_k}} - S_{n_k}| > \frac{\eta}{\sqrt{2c}}, \sigma\sqrt{2 \cdot 2ca_{n_k}d_{n_k}}\bigg)$$

$$+ 2P(|S_{n_k+a_{n_k}} - S_{n_k}| > (1+\eta)\sigma\sqrt{2a_{n_k}d_{n_k}}).$$

Summing the three probabilities over $k$ and recalling (3.24) tells us that the total sum converges whenever

$$\min\bigg\{\frac{\eta}{\sqrt{2c(1+\eta)^3}}, \frac{\eta}{\sqrt{2c}}, 1+\eta\bigg\} > 1.$$

Since we can choose $c > 0$ arbitrarily small, we finally conclude that for any $\eta > 0$, we have

$$\sum_k P\bigg(\max_{n_k \leq n \leq n_{k+1}} \frac{S_{n+a_n} - S_n}{\sqrt{2a_n d_n}} > (1+6\eta)\sigma\bigg) < \infty,$$

implying the upper inequality for the entire sequence, that is,

$$\limsup_{n\to\infty} \frac{T_{n,n+n/L(n)}}{\sqrt{2a_n d_n}} \leq \sigma \qquad \text{a.s.}$$

### 3.8. Necessity

By the zero-one law, the probability that the lim sup is finite is 0 or 1, hence, being positive, it equals 1. Consequently, (see [9], page 438),

$$\limsup_{n\to\infty} \frac{|X_n|}{\sqrt{a_n d_n}} < \infty \qquad \text{a.s.},$$

from which, in the case $a_n d_n \leq f(n)$, it follows, via the second Borel–Cantelli lemma and the i.i.d. assumption, that

$$\sum_{n=1}^{\infty} P(|X_n| > \sqrt{f(n)}) = \sum_{n=1}^{\infty} P(X^2 > f(n)),$$

which verifies (2.1).

If (2.1) is weaker than finite variance, that is, if $a_n d_n \geq f(n)$, then, by following Feller's proof [4] (see also, for example, [5], Section 8.4) of the necessity in the LIL with only obvious changes involving the replacing of sums by windows, we may conclude that (2.1) – that is, finite variance – also holds in this case.

An application of the sufficiency part then tells us that (2.1) holds with $\sigma^2 = \text{Var}\,X$.

Finally, if $EX = \mu$, then, by the law of large numbers, $S_n/f(n) \sim \mu n/f(n) \to \mu \cdot \infty$ as $n \to \infty$, which forces $\mu$ to be equal to zero.

## 4. Proofs of the Corollaries in Section 2

### 4.1. Proof of Corollary 2.1

In this case, $a_n = n/(\log n)^p$ (for $n \geq 9$) so that

$$d_n = \log \frac{n}{n/(\log n)^p} + \log \log n = (p+1)\log \log n,$$

that is, $f(n) = (p+1)n \log \log n/(\log n)^p$. It follows that $f^{-1}(n) \sim \frac{n(\log n)^p}{(p+1)\log \log n}$ as $n \to \infty$ so that (2.1) turns out as

$$EX^2 \frac{(\log^+ |X|)^p}{\log^+ \log^+ |X|} < \infty.$$

It remains to verify that $xL'(x)/L(x)$ is decreasing. Now, $L(x) = (\log x)^p$ and $L'(x) = x^{-1}p(\log x)^{p-1}$ so that $xL'(x)/L(x) = p(\log x)^{-1}$, which indeed decreases.

### 4.2. Proof of Corollary 2.2

Thus, $a_n = n/\log_m(n)$ for $n$ sufficiently large, so that

$$d_n = \log \log_m(n) + \log \log n = \log_{m+1}(n) + \log \log n \sim \log \log n \qquad \text{as } n \to \infty.$$

Since $a_n d_n > n$, we have $f^{-1}(n) = n$ as $n \to \infty$, which implies that finite variance is the appropriate necessary assumption.

As for (1.2), this time, $L(x) = \log_m x$ and $L'(x) = x^{-1} \prod_{i=1}^{m-1}(\log_i x)^{-1}$ so that $xL'(x)/L(x) = \prod_{i=1}^{m}(\log_i x)^{-1}$, which indeed decreases.

## 5. Further examples

In this section we provide some additional examples to illustrate Theorem 2.1. As in Section 2 we omit stating converse results.

The first example mixes powers of logarithms and iterated logarithms.

*Example 5.1.* Let, for $n \geq 9$, $a_n = n(\log \log n)^q/(\log n)^p$, $p, q > 0$, which means that the slowly varying function part is $L(n) = (\log n)^p/(\log \log n)^q$. Differentiation and some algebraic simplification yield that $xL'(x)/L(x) = p/\log x - q/(\log x \log_2 x)$, which is ultimately decreasing. Moreover,

$$\begin{aligned}d_n &= \log\left(\frac{n(\log \log n)^q}{n/(\log n)^p}\right) + \log \log n \\ &= (p+1)\log \log n - q \log \log \log n \sim (p+1)\log \log n \qquad \text{as } n \to \infty\end{aligned}$$

so that $f(n) = (p+1)n(\log\log n)^{q+1}/(\log n)^p$, which implies that $f^{-1}(n) \sim Cn(\log n)^p/(\log\log n)^{q+1}$ as $n \to \infty$. The following result emerges.

**Corollary 5.1.** *Suppose that $X, X_1, X_2, \ldots$ are i.i.d. random variables with mean 0 and finite variance $\sigma^2$, and let $T_{n,n+k} = \sum_{j=n+1}^{n+k} X_j$. If, for some $p, q > 0$,*

$$EX^2 \frac{(\log^+ |X|)^p}{(\log^+ \log^+ |X|)^{q+1}} < \infty, \tag{5.1}$$

*then*

$$\limsup_{n \to \infty} \frac{T_{n,n+n(\log\log n)^q/(\log n)^p}}{\sqrt{2(p+1)(n/(\log n)^p)(\log\log n)^{q+1}}} = \sigma \quad \text{a.s.} \tag{5.2}$$

The previous conclusion also holds, in fact, for $q = 0$, in which case Corollary 5.1 reduces to Corollary 2.1 since the $\log\log$-contribution is of a lower order of magnitude. However, the case $p = 0$ requires a separate treatment.

*Example 5.2.* Let, for $n \geq 9$, $a_n = n/(\log\log n)^q$, $q > 1$. Now, $L(x) = (\log_2 x)^q$ gives $xL'(x)/L(x) = q/(\log x \log_2 x)$, which is decreasing. Moreover,

$$d_n = \log\left(\frac{n}{n/(\log\log n)^q}\right) + \log\log n = q\log\log\log n + \log\log n \sim \log\log n \quad \text{as } n \to \infty,$$

that is, $f(n) = n(\log\log n)^{1-q}$, and $f^{-1}(n) \sim n(\log\log n)^{q-1}$ as $n \to \infty$.

**Corollary 5.2.** *Suppose that $X, X_1, X_2, \ldots$ are i.i.d. random variables with mean 0 and finite variance $\sigma^2$, and let $T_{n,n+k} = \sum_{j=n+1}^{n+k} X_j$. If, for some $q > 1$,*

$$EX^2(\log^+ \log^+ |X|)^{q-1} < \infty, \tag{5.3}$$

*then*

$$\limsup_{n \to \infty} \frac{T_{n,n+n/(\log\log n)^q}}{\sqrt{2n(\log\log n)^{1-q}}} = \sigma \quad \text{a.s.} \tag{5.4}$$

*Example 5.3.* Let $a_n = n/\exp\{\sqrt{\log n}\}$, $n \geq 1$. Since $L(x) = \exp\{\sqrt{\log x}\}$, we have $xL'(x)/L(x) = (\log x)^{-1/2}/2$, which is decreasing. Moreover,

$$d_n = \log\exp\{\sqrt{\log n}\} + \log\log n = \sqrt{\log n} + \log\log n \sim \sqrt{\log n} \quad \text{as } n \to \infty,$$

which gives $f(n) \sim n\sqrt{\log n}/\exp\{\sqrt{\log n}\}$ as $n \to \infty$ so that

$$f^{-1}(n) \sim n\exp\{\sqrt{\log n + 1/2}\}/\sqrt{\log n} \quad \text{as } n \to \infty.$$

The following conclusion therefore holds.

**Corollary 5.3.** *Suppose that* $X, X_1, X_2, \ldots$ *are i.i.d. random variables with mean 0 and finite variance* $\sigma^2$, *and let* $T_{n,n+k} = \sum_{j=n+1}^{n+k} X_j$. *If*

$$EX^2 \frac{\exp\{\sqrt{2\log^+ |X|}\}}{\sqrt{\log^+ |X|}} < \infty, \tag{5.5}$$

*then*

$$\limsup_{n\to\infty} \frac{T_{n,n+n/\exp\{\sqrt{\log n}\}}}{\sqrt{2(n/\exp\{\sqrt{\log n}\})\sqrt{\log n}}} = \sigma \quad a.s. \tag{5.6}$$

The following example is a more general version.

**Example 5.4.** Let, for $n \geq 1$, $a_n = n(\log n)^\gamma / \exp\{(\log n)^\beta\}$, where $0 < \beta < 1$ and $\gamma \in \mathbb{R}$. Thus, $L(x) = (\log x)^{-\gamma} \exp\{(\log x)^\beta\}$ and $xL'(x)/L(x) = \beta(\log x)^{\beta-1} - \gamma/\log x$, which is ultimately decreasing. Furthermore,

$$d_n = \log \exp\{(\log n)^\beta\} - \log\log((\log n)^\gamma) + \log\log n \sim (\log n)^\beta \qquad \text{as } n \to \infty,$$

which gives

$$f(n) \sim n(\log n)^{\beta+\gamma} / \exp\{(\log n)^\beta\} \qquad \text{as } n \to \infty.$$

It follows that for $0 < \beta < 1/2$,

$$f^{-1}(n) \sim n \frac{\exp\{(\log n)^\beta\}}{(\log n)^{\beta+\gamma}} \qquad \text{as } n \to \infty,$$

for $1/2 \leq \beta < 2/3$,

$$f^{-1}(n) \sim n \frac{\exp\{(\log n)^\beta + \beta(\log x)^{2\beta-1}\}}{(\log n)^{\beta+\gamma}} \qquad \text{as } n \to \infty,$$

and so on. The following conclusion holds.

**Corollary 5.4.** *Suppose that* $X, X_1, X_2, \ldots$ *are i.i.d. random variables with mean 0 and finite variance* $\sigma^2$, *and let* $T_{n,n+k} = \sum_{j=n+1}^{n+k} X_j$. *If* $0 < \beta < 1/2$ *and*

$$EX^2 \frac{\exp\{(2\log^+ |X|)^\beta\}}{(\log^+ |X|)^{\beta+\gamma}} < \infty, \tag{5.7}$$

*then*

$$\limsup_{n\to\infty} \frac{T_{n,n+n(\log n)^\gamma/\exp\{(\log n)^\beta\}}}{\sqrt{2n(\log n)^{\gamma+\beta}/\exp\{(\log n)^\beta\}}} = \sigma \quad a.s. \tag{5.8}$$

## 6. A simplified Theorem 2.1

As mentioned in the Introduction, this section concerns a weaker result, the proof of which is *much* easier, in that only one truncation is made. However, minimal moment conditions are not obtained.

**Theorem 6.1.** *Suppose that $X, X_1, X_2, \ldots$ are i.i.d. random variables with mean 0 and finite variance $\sigma^2$, and let $T_{n,n+k} = \sum_{j=n+1}^{n+k} X_j$. Define sequences $\{a_n\}$ and $\{d_n\}$ as in Section 2 and set*

$$b_n = \sqrt{\frac{a_n}{d_n}}.$$

*If $EX = 0$ and*

$$Eb^{-1}(|X|) < \infty, \tag{6.1}$$

*then*

$$\limsup_{n \to \infty} \frac{T_{n,n+a_n}}{\sqrt{2a_n d_n}} = \sigma \quad a.s. \tag{6.2}$$

### 6.1. Proof of Theorem 6.1

Let

$$X_n'' = X_n I\left\{|X_n| > \frac{\sigma\delta}{\varepsilon} b_n\right\} = X_n - X_n'.$$

*The contribution of $T_{n,n+a_n}'$*

This part requires no change. In other words, we first let $\{n_k, k \geq 1\}$, where $n_k \nearrow \infty$ as $k \to \infty$, satisfy (3.10), and apply (3.4) to $\{X_{n_k}, k \geq 1\}$ to obtain

$$\sum_{k=1}^{\infty} P(|T_{n_k,n_k+a_{n_k}}'| > \varepsilon\sqrt{2a_{n_k} d_{n_k}}) < \infty$$

for the convergence part.

After this, we choose a sparser subsequence $\{n_k, k \geq 1\}$, where $n_k \nearrow \infty$ as $k \to \infty$, satisfying (3.12) and apply (3.6) to obtain

$$\sum_{k=1}^{\infty} P(T_{n_k,n_k+a_{n_k}}' > \varepsilon\sqrt{2a_{n_k} d_{n_k}}) = \infty$$

for the divergence part.

*The contribution of $T''_{n,n+a_n}$*

Since the $X''$'s have changed, we can use the stronger LSL argument here (see [6, 7]). Namely, in order for the $|T''_{n,n+a_n}|$'s to surpass the level $\eta\sqrt{a_n d_n}$ infinitely often, it is necessary that infinitely many of the $X''$'s are non-zero. However, the latter event has zero probability in view of the first Borel–Cantelli lemma since

$$\sum_{n=1}^{\infty} P\left(|X_n| > \frac{\sigma\delta}{\varepsilon}b_n\right) = \sum_{n=1}^{\infty} P\left(|X| > \frac{\sigma\delta}{\varepsilon}b_n\right) < \infty$$

if and only if (6.1) holds.

*Completing the proof*

From this point on, the arguments are identical to those of the proof of Theorem 2.1. We therefore omit the details.

### 6.2. Revisiting some examples

**Corollary 2.1 revisited for $p=1$.** With $a_n = n/\log n$ and $b_n = \sqrt{\frac{n/\log n}{\log\log n}}$, the conclusion from Theorem 6.1 is, as before,

$$\frac{T_{n,n+n/\log n}}{\sqrt{4(n/\log n)\log\log n}} = \sigma \qquad a.s.,$$

*however, provided (6.1) holds, that is, provided*

$$E(X^2 \log^+ |X| \log^+ \log^+ |X|) < \infty,$$

*since $b^{-1}(n) \sim n^2 \log n \log\log n$. This should be compared with (2.4),*

$$EX^2 \frac{\log^+ |X|}{\log^+ \log^+ |X|} < \infty.$$

**Corollary 2.2 revisited for $m=2$.** With $m=2$, $a_n = n/\log\log n$ and $b_n = \sqrt{\frac{n/\log\log n}{\log\log n}}$, *the conclusion from Theorem 6.1 is, as before,*

$$\frac{T_{n,n+n/\log\log n}}{\sqrt{2n}} = \sigma \qquad a.s.,$$

*however, provided (6.1) holds, that is, provided*

$$EX^2(\log^+ \log^+ |X|)^2 < \infty,$$

*which should be compared with the optimal one which was finite variance.*

**Corollary 5.1 revisited.** *Here, $a_n = n(\log\log n)^q/(\log n)^p$, $p, q > 0$. With*

$$b_n = \sqrt{\frac{n(\log\log n)^q/(\log n)^p}{\log\log n}} = \sqrt{\frac{n(\log\log n)^{q-1}}{(\log n)^p}},$$

*assumption (6.1) turns out as*

$$EX^2 \frac{(\log^+ |X|)^p}{(\log^+ \log^+ |X|)^{q-1}} < \infty,$$

*to be compared with the weaker*

$$EX^2 \frac{(\log^+ |X|)^p}{(\log^+ \log^+ |X|)^{q+1}} < \infty.$$

## Acknowledgement

We wish to thank an anonymous referee of [7] for drawing our attention to the problem of analyzing the limiting behavior of windows of the form (1.1).